\documentclass[12pt]{article}
\usepackage{amsmath}
\begin{document}
\title{An illustration of a paradox about the idoneal, or suitable,
numbers\footnote{Delivered to the St.--Petersburg Academy April 20, 1778.
Originally published as
\emph{
Illustratio paradoxi circa progressionem numerorum idoneorum sive congruorum},
Nova Acta Academiae Scientarum Imperialis Petropolitinae \textbf{15} (1806),
29--32, and
republished in \emph{Leonhard Euler, Opera Omnia}, Series 1:
Opera mathematica,
Volume 4, Birkh\"auser, 1992. A copy of the original text is available
electronically at the Euler Archive, at http://www.eulerarchive.org. This paper
is E725 in the Enestr\"om index.}}
\author{Leonhard Euler\footnote{Date of translation: July 17, 2005.
Translated from the Latin
by Jordan Bell, 3rd year undergraduate in Honours Mathematics, School of Mathematics and Statistics, Carleton University,
Ottawa, Ontario, Canada.
Email: jbell3@connect.carleton.ca.
This translation was written
during an NSERC USRA supervised by Dr. B. Stevens.
}}
\date{}
\maketitle

I. An outstanding paradox stands upon this,
for although the idoneal numbers are shaped and proceed according
to a certain law, 
the multitude of which however are not infinite yet are extended 
even onto 65 terms, concerning this paradox
I have recorded so far no more
of this type in the succession which has been observed;
yet neither on the other hand
has it been permitted by me to make firm
a finite number of terms, except that after the 65th term,
which is 1848, none thereafter have been bestowed, even though I have
continued the examination up to 10000 and beyond.

II. Neither still is any other route seen to stand open to
this outstanding paradox which is to be demonstrated. 
On account of this, it will be not insufficiently conveyed of
the light
for this most recondite things, 
since at least for a certain type of these numbers, namely
the squares, it will be possible to be demonstrated for
the multide of these forms to be terminated, that is
for it to not be possible for other square numbers to occur
in the succession of idoneal numbers, aside from the 
first five $1,4,9,16$ and 25,
insofar as I will demonstrate this in the following
way by a certain rule for this progression.

III. Consequently we transfer the rule for finding idoneal numbers,
having been related in a certain published paper, specifically
for square numbers, which will thus be enunciated in the following way:

From the sequence of all the square numbers,
for each prime number $p$, all the numbers contained
in the form $px \quad yy$ [sic] and larger than
$\frac{1}{4}pp$ are excluded, aside from those as $pp-yy$,
since if this were done for each prime number, those that will be
left behind from the sequence of square numbers are the very same
as those which are idoneal. Indeed among the primes, 
we have seen in place of two the square 4 of it ought to be taken.
With then the form $px-yy$ involving no square numbers,
then nothing which has been removed should have a place.

IV. Similarly for the prime number $p=3$ it turns out that
indeed the form $3x-yy$ envelops no square numbers, 
insofar as this likewise would hold for all numbers of the form
$p=4n-1$. For if the form $(4n-1)x-yy$ were a square,
valued at $zz$, the sum of the two squares $yy+zz$ would be divisible
by $4n-1$, which has been noted to be impossible; from this
it is realized for our $p$ not to be left as other prime numbers unless
they are contained in the form $4n+1$.

V. Therefore it will be $p=5$, such that from the sequence of squares,
those which are contained in the form $5x-yy$ and which
surpass $\frac{1}{4}pp=6\frac{1}{4}$ should be excluded, excepting
however those contained in the form $25-yy$, which are 9 and 16,
from which all larger squares contained in the form
$5x-yy$ will ought to be removed. Therefore with it that
all squares not divisible by 5 are either of
the form $5x-1$ or $5x-4$, it is now evident for all squares
not divisible by 5 and simultaneously greater than $6 \frac{1}{4}$
to ought to be excluded from the sequence of
all the square numbers; then with this having been done, the following
sequence of squares will be left behind: $1,4,9,16,25,10^2,15^2,20^2$.
 Namely here after 16 no other squares are left, unless the roots
of them are divisible by 5.

VI. The following number of the form $4x+1$ is $p=13$, from
which those numbers contained in the form $13x-yy$ should be excluded
which are larger than $42 \frac{1}{4}$, excepting however those
which are contained in the form $169-yy$, which are 25 and 144.
Therefore those square numbers which are to be
excluded, greater than $42 \frac{1}{4}$, are contained
in the form $13x-yy$; this form clearly
contains all squares whose roots are not divisible by 13,
and by this exclusion, no others remain after 25 unless divisible by 13.
Moreover, by the preceding
condition no others are left unless divisible by 5; from this
therefore if all those not divisible by 13 are taken away
which are less than $42 \frac{1}{4}$, aside from 25, no other
squares  are left behind, unless they are simultaneously divisible by 5 and 13,
from which therefore all are contained in the form $(65 \alpha)^2$;
therefore the surviving square numbers are: $1,4,9,16,25,65^2,130^2,
195^2,260^2,$ etc.

VII. The following prime number of the form $4n+1$ is $p=17$,
and hence the form of those which are to be excluded will
be $17x \quad yy$ [sic], so far as they are larger than
$\frac{1}{4}pp=72 \frac{1}{4}$, excepting however those
which are contained in the form $17^2-yy$, which are
$15^2$ and $8^2$,
which indeed by the preceding conditions have been deleted.
Therefore, all squares not divisible by 17 should be deleted;
from this it is apparent for no others after 25 to be left, unless
they at once are divisible by 5, 13 and 17, which are therefore
contained in the form $(5,13,17)^2$,
of which the first is $1105^2$.

VIII. The following prime number of the form $4n+1$ is $p=29$,
from which the form containing the numbers which are to be 
exluded is $29x-yy$, namely in that it contains squares, such that
it would be $29x=yy+zz$. Indeed, clearly
this form contains all squares not divisible by 29;
thus with the previous ones which have been eliminated none remain,
unless in they are contained in the form $(5,13,17,29\alpha)^2$, of
which the smallest is 32045.

IX. So if in this way we evaluate the successive prime numbers of
the form $4n+1$, it is evident that after the initial five squares
$1,4,9,16$ and 25 for no other idoneal to occur on to infinity.

X. Therefore if the question is instituted about the square numbers,
which are simultaneously idoneals, it has now been firmly
demonstrated for such numbers not to be so ascribed,
except for these five: $1,4,9,16,25$; 
from this it is now permitted to understand clearly, to the
extent that, not hindering the law of this progression,
a great number of all the idoneal numbers are plainly able to be
terminated, and perhaps sometime the others will be able to be
demonstrated in a similar way.

XI. In this demonstration we have assumed for all squares to be contained
in the form $px-yy$, which are not divisible by the number $p$.
For $p=4n+1$ is always the sum of two squares, which would be
$aa+bb$, so that we will thus have
$(aa+bb)x-yy=zz$, from which, by taking $x=ff+gg$,
it is gathered, in the form $px-yy$ for the number
$y$ to be bound to be prime to $p$.

\end{document}